\input amstex\documentstyle{amsppt}  
\pagewidth{12.5cm}\pageheight{19cm}\magnification\magstep1
\topmatter
\title{Characters of simplylaced nonconnected groups versus characters of nonsimplylaced connected groups}
\endtitle
\author{Shrawan Kumar, George Lusztig and Dipendra Prasad}\endauthor
\abstract{Let $G$ be a connected, simply-connected, almost simple semisimple group over $\bold C$ of simplylaced 
type and let $\sigma$ be a nontrivial diagram automorphism of $G$.  Let $G\langle\sigma\rangle$ be the 
(disconnected) group generated by $G$ and $\sigma$. As a consequence of a theorem of Jantzen the character of an 
irreducible representation of $G\langle\sigma\rangle$ (which remains 
 irreducible on $G$) on $G\sigma$ can be expressed in 
terms of a character of an irreducible representation of a certain connected 
simply-connected semisimple group 
$G_{\sigma}$ of nonsimplylaced type. We show how Jantzen's theorem can be deduced from properties of the 
canonical bases.}\endabstract
\endtopmatter   
\document

\define\frl{\forall}

\define\lb{\linebreak}

\define\part{\partial}

\define\ra{\rangle}

\define\m{\mapsto}

\define\la{\langle}

\define\T{\times}

\define\nl{\newline}
\redefine\i{^{-1}}

\define\ot{\otimes}

\define\Hom{\text{\rm Hom}}

\define\tr{\text{\rm tr}}

\define\a{\alpha}

\define\et{\eta}

\define\s{\sigma}

\redefine\l{\lambda}

\define\vp{\varpi}

\define\BB{\bold B}
\define\CC{\bold C}

\define\NN{\bold N}

\define\ZZ{\bold Z}

\define\co{\Cal O}

\define\sh{\sharp}

\define\che{\check}
\define\cha{\che{\a}}

\define\FSS{FSS}
\define\FRS{FRS}
\define\JA{Ja}
\define\QG{L1}
\define\CL{L2}
\define\MO{Mo}
\define\NA{N1}
\define\NAA{N2}
\define\NS{NS}

Let $G$ be a connected, simply-connected, almost simple algebraic group of simplylaced type over $\CC$.
Let $T$ be a maximal torus of $G$. Let $x_i:\CC@>>>G$, $y_i:\CC@>>>G$ $(i\in I)$ be homomorphisms which together 
with $T$ form a pinning (\'epinglage) of $G$. We fix a nontrivial automorphism $\s$ of $G$ such that $\s(T)=T$, 
and such that for some permutation $i\m\tilde{i}$ of $I$ we have $\s(x_i(a))=x_{\tilde{i}}(a)$, 
$\s(y_i(a))=y_{\tilde{i}}(a)$ for all $a\in\CC$. For $i\in I$ we write $\s(i)=\tilde{i}$. Let $\la\s\ra$ be the 
finite subgroup of the automorphism group of $G$ generated by $\s$ and let $G\la\s\ra$ be the semidirect product 
of $G$ with $\la\s\ra$.

Let $X$ be the group of characters $T@>>>\CC^*$; let $Y$ be the group of one parameter subgroups $\CC^*@>>>T$ and
let $\la,\ra:Y\T X@>>>\ZZ$ be the standard pairing. For $i\in I$ we define $\a_i\in X$ by
$x_i(\a_i(t))=tx_i(1)t\i$, $y_i(\a_i(t)\i)=ty_i(1)t\i$ for all $t\in T$. This is a root of $G$. Let $\cha_i\in Y$
be the corresponding coroot. Note that 

(a) $(Y,X,\la,\ra,\cha_i,\a_i (i\in I))$ 
\nl
is the root datum of $G$.
Now $\s$ induces automorphisms of $X,Y$ denoted again by $\s$; these are compatible with $\la,\ra$ and we have
$\s(\a_i)=\a_{\s(i)},\s(\cha_i)=\cha_{\s(i)}$ for $i\in I$. Let 
$X^+=\{\l\in X;\la\cha_i,\l\ra\in\NN\frl i\in I\}$.

We set $Y_\s=Y/(\s-1)Y$, ${}^\s X=\{\l\in X;\s(\l)=\l\}$. Note that $\la,\ra:Y\T X@>>>\ZZ$ induces a perfect 
pairing $Y_\s\T{}^\s X@>>>\ZZ$ denoted again by $\la,\ra$. Let $I_\s$ be the set of $\s$-orbits on $I$. For any 
$\co\in I_\s$ let $\cha_\co\in Y_\s$ be the image of $\cha_i$ under $Y@>>>Y_\s$ where $i$ is any element of 
$\co$. Since $\{\cha_i;i\in I\}$ is a $\ZZ$-basis of $Y$ we see that $\{\cha_\co;\co\in I_\s\}$ is a $\ZZ$-basis 
of $Y_\s$. For any $\co\in I_\s$ let $\a_\co=2^h\sum_{i\in\co}\a_i\in{}^\s X$ where $h$ is the number of unordered
pairs $(i,j)$ such that $i,j\in\co$, and $\a_i+\a_j$ is a root. Note that $h=0$ except when $G$ is of type 
$A_{2n}$ when $h=0$ for all $\co$ but one and $h=1$ for one $\co$. Note that

(b) $(Y_\s,{}^\s X,\la,\ra,\cha_\co,\a_\co (\co\in I_\s))$ 
\nl
is a root datum, see \cite{\JA, p.29}. Let 
${}^\s X^+=\{\l\in{}^\s X;\la\cha_\co,\l\ra\in\NN\frl\co\in I_\s\}={}^\s X\cap X^+$. Let $G_\s$ be the connected 
semisimple group over $\CC$ with root datum (b). By definition, $G_\s$ is provided with an \'epinglage 
$(T_\s,x_\co,y_\co\quad(\co\in I_\s))$ where $T_\s:=\CC^*\ot Y_\s=T/\{\s(t)t\i;t\in T\}$ is a maximal torus of 
$G_\s$ and $x_\co:\CC@>>>G_\s$, $y_\co:\CC@>>>G_\s$ satisfy
$x_\co(\a_\co(t_1))=t_1x_\co(1)t_1\i$, $y_\co(\a_\co(t_1)\i)=t_1y_\co(1)t_1\i$ for all $t_1\in T_\s$. 
(We have ${}^\s X=\Hom(T_\s,\CC^*)$ canonically.) 

Note that $G_\s$ is simply-connected and that $G_\s\cong {}^L(({}^\s ({}^LG))^0)$ where ${}^L()$ denotes the 
Langlands dual group and ${}^\s ({}^LG)$ denotes the fixed point set of the automorphism of ${}^LG$ induced by $\s$. 
Now $G_\s$ is only isogenous to ${}^L({}^\s G)$ where ${}^\s G$ is the fixed point set of $\s:G@>>>G$.

Let $\l\in{}^\s X^+$. We can view $\l$ both as a character of $T$ and as a character of $T_\s$. Let $V$ (resp. 
$V'$) be a finite dimensional complex irreducible representation of $G$ (resp. $G_\s$) with a non-zero vector 
$\et$ (resp. $\et'$) such that $x_i(a)\et=0$ for all $i\in I,a\in\CC$ (resp. $x_\co(a)\et'=0$ for all 
$\co\in I_\s,a\in\CC$) and $t\et=\l(t)\et$ for all $t\in T$ (resp. 
$t'\et'=\l(t')\et'$ for all $t'\in T_\s$). Now 
$V$ can be regarded as a representation of $G\la\s\ra$ whose restriction to $G$ is as above and on which the 
action of $\s$ satisfies $\s(\et)=\et$. 

Let $\mu\in X$. Let $V_\mu=\{x\in V;tx=\mu(t)x\quad \frl t\in T\}$. Note that $\s:V@>>>V$ permutes the weight 
spaces $V_\mu$ among themselves. A weight space $V_\mu$ is $\s$-stable if and only if $\mu\in{}^\s X$; in this 
case $\mu$ can be viewed as a character of $T_\s$  and we set $V'_\mu=\{x'\in V';t'x'=\mu(t')x'\frl t'\in T_\s\}$.

\proclaim{Theorem (Jantzen \cite{\JA, Satz 9})} For $\mu\in {}^\s X, \tr(\s:V_\mu@>>>V_\mu)=\dim V'_\mu$.
\endproclaim
\proclaim{Corollary} Let $\vp:T@>>>T_\s$ be the canonical homomorphism. For any $t\in T$ we have 
$\tr(t\s:V@>>>V)=\tr(\vp(t),V')$.
\endproclaim
The corollary describes completely the character of $V$ on $G\s$ in terms of the character of $V'$ since any 
semisimple element in $G\s$ is $G$-conjugate to an element of the form $t\s$ with $t\in T$. Note also that there 
is a well defined bijection between the set of semisimple $G$-conjugacy classes in $G\s$ and the set of
semisimple $G_\s$-conjugacy classes in $G_\s$ which for any $t\in T$ maps the $G$-conjugacy class of $t\s$ to the
$G_\s$-conjugacy class of $\vp(t)$; see \cite{L2, 6.26}, \cite{\MO}.

We now show (assuming that $G$ is not of type $A_{2n}$) how Jantzen's theorem can be deduced from properties of 
canonical bases in \cite{\QG}. According to \cite{\QG}, $V$ has a canonical basis $B_\l$ and $V'$ has a canonical
basis $B'_\l$. Also, $B_\l$ (resp. $B'_\l$) can be naturally viewed as a subset of $\BB$ (resp. $\BB'$), the 
canonical basis of the $+$ part of the universal enveloping algebra attached to the root datum (a) (resp. (b)). 
Now $\s$ acts naturally on $\BB$ (preserving the subset $B_\l$) and \cite{\QG, Theorem 14.4.9} provides a
canonical bijection between $\BB'$ and the fixed point set of $\s$ on $\BB$. (This theorem is applicable since the
Cartan datum of (b) is obtained from the Cartan datum of (a) by the general ``folding" procedure \cite{\QG, 14.1} 
which appplies to any simplylaced Cartan datum of not necessarily finite type together with an admissible 
automorphism; here we use that $G$ is not of type $A_{2n}$.) This restricts to a bijection between $B'_\l$ and the
fixed point set ${}^\s B_\l$ of $\s$ on $B_\l$. Next we note that $B_\l$ (resp. $B'_\l$) is compatible with the 
decomposition of $V$ (resp. $V'$) into weight spaces and from the definitions we see that the bijection above 
carries $B'_\l\cap V'_\mu$ bijectively onto ${}^\s B_\l\cap V_\mu$. Since $B_\l\cap V_\mu$ is a basis of $V_\mu$ 
which is $\s$-stable we have $\tr(\s:V_\mu@>>>V_\mu)=\sh({}^\s B_\l\cap V_\mu)$. Using the bijection above this 
equals $\sh(B'_\l\cap V'_\mu)$ and this is equal to $\dim V'_\mu$ since $B'_\l\cap V'_\mu$ is a basis of $V'_\mu$.
This gives the desired result.

We refer the reader to \cite{\FSS, \FRS, \NS, \NA, \NAA} for other approaches to Jantzen's theorem. We thank 
S. Naito for pointing out these references to us. 

The first two authors were supported in part by the National Science Foundation. The third author thanks the 
Institute for Advanced Study where this work was done, and gratefully acknowledges receiving support through 
grants to the Institute by the Friends of the Institute, and the von Neumann Fund.

Addresses:

S.K.: Department of Mathematics, University of North Carolina, Chapel Hill, NC 27599-3250, USA

G.L.: Department of Mathematics, M.I.T., Cambridge, MA 02139, USA

D.P.: School of Mathematics, Tata Institute of Fundamental Research, Colaba,
Mumbai 400005, India, and The Institute for Advanced Study, Princeton, NJ 08540, USA
\enddocument